\documentclass{article}
\usepackage[utf8]{inputenc}

\usepackage{amsmath,amssymb}
\usepackage{amsthm}
\usepackage{enumitem}
\usepackage{multicol}
\usepackage{hyperref}
\usepackage{array}
\usepackage{url}
\usepackage{listings}
\usepackage{color}
\usepackage{tikz,tkz-graph}
\usepackage{mwe}
\usepackage{booktabs}
\usepackage{colortbl}
\usepackage{array}
\tikzset{
every node/.style={circle, draw, inner sep=2pt},
every picture/.style={thick}
}
\usepackage{graphicx}
\usepackage{subcaption} 
\usetikzlibrary{matrix,arrows,calc}
\usetikzlibrary{decorations.pathreplacing,decorations.markings}

\hypersetup{
 colorlinks=true,
 linkcolor=blue,
 filecolor=blue,
 citecolor = blue,      
 urlcolor=blue,
 }

\theoremstyle{definition}

\title{Evolutive sandpiles}
\author{Carlos A. Alfaro$^a$, Juan Pablo Serrano$^b$ and Ralihe R. Villagrán$^c$
\\ \\
{\small $^a$Banco de M\'exico} \\
{\small Mexico City, Mexico} \\
{\small {\tt alfaromontufar@gmail.com,carlos.alfaro@banxico.org.mx}} \\
\\
{\small $^b$Departamento de Matem\'aticas}\\
{\small Centro de Investigaci\'on y de Estudios Avanzados del IPN}\\
{\small Apartado Postal 14--740 }\\
{\small 07000 Mexico City, Mexico. } \\
{\small {\tt jpserranop@math.cinvestav.mx}} \\
\\
{\small $^c$Department of Mathematical Sciences,}\\ 
{\small Worcester Polytechnic Institute, Worcester, USA}\\
{\small {\tt rvillagran@wpi.edu}}\\
}
\date{}

\begin{document}

\maketitle

\begin{abstract}

The Abelian sandpile model was the first example of a self-organized critical system studied by Bak, Tang and Wiesenfeld.
The dynamics of the sandpiles occur when the grains topple over a graph.
In this study, we allow the graph to evolve over time and change the topology at each stage.
This turns out in the occurrence of phenomena impossible in the classical sandpile models.
For instance, configurations over evolutive graphs that are always unstable.
We also experiment with the stabilization of configurations with a large number of grains at the center over evolutive graphs, this allows us to obtain interesting fractals.
Finally, we obtain power laws associated with some evolutive graphs.
\end{abstract}

\emph{Keywords:} Evolutive sandpiles, fractals, power laws, self-organized critical systems

\section{Introduction}
The Abelian sandpile model was the first example of a self-organized critical system studied by Bak, Tang and Wiesenfeld.
It was the first example of a {\it self-organized critical system}, which attempts to explain the occurrence of power laws in many natural phenomena \cite{Bak96} ranging on different fields like geophysics \cite{geoph}, optimization \cite{optimi1,optimi2,ALFARO2023113356}, economics \cite{Biggs99} and neuroscience \cite{brain, WILTING2019105}.
An exposure to self-organized-critically is provided in \cite{Bak96}, and a more recent re-account can be found in \cite{Watkins25years}. For a short and clear explanation of the sandpile model we can refer the reader to \cite{WhatIs} for instance.

The dynamics of the {\it Abelian sandpile model}, which was first studied by Bak, Tang and Wiesenfeld in \cite{PhysRevLett.59.381}, is carried out on a connected graph $G=(V,E)$ with a special vertex $q\in V$, called {\it sink}.
Let $\mathbb{N}$ denote the set of non-negative integers and let $\mathbb{Z}$ denote the integer numbers.
In the sandpile model, a configuration on $(G,q)$ is a vector ${ c}\in \mathbb{N}^{V}$, in which the entry ${c}_v$ is associated with the number of {\it grains of sand} (depending on the context it can also be used {\it chips} or {\it dollars} instead) placed on vertex $v$.
The sink vertex is used to collect the sand getting out of the system.
Let $\widetilde{V}$ denote the set of non-sink vertices.
Two configurations $c$ and $d$ are {\it equal} if $c_v=d_v$ for each non-sink vertex $v\in \widetilde{V}$.
A non-sink vertex $v$ is called \textit{stable} if ${c}_v$  is less than its {\it degree} $d_G(v)$, and {\it unstable}, otherwise.
Thus, a configuration is called \textit{stable} if every non-sink vertex is stable.
The {\it toppling rule} in the dynamics of the model consists of selecting an unstable non-sink vertex $u$ and moving $d_G(u)$ grains of sand from $u$ to its neighbors, in which each neighbor $v$ receives $m_{(u,v)}$ grains of sand, where $m_{(u,v)}$ denote the number of edges between $u$ and $v$.
Note toppling vertex $v_i$ in configuration ${c}$ corresponds to the subtraction of the $i$-{\it th} row of the Laplacian matrix to ${c}$. 
Recall the Laplacian matrix $L(G)$ of a graph $G$ is such that the $(u,v)$-entry of $L(G)$ is defined as 
\[
L(G)_{u,v}=
\begin{cases}
\deg_G(u) & \text{if } u=v,\\
-m(u,v) & \text{otherwise.}
\end{cases}
\]
Starting with any unstable configuration and toppling unstable vertices repeatedly, we will always obtain \cite[Theorem 2.2.2]{MR3889995} a stable and unique configuration after a finite sequence of topplings.
The stable configuration obtained from the configuration ${ c}$ will be denoted by $s({c})$.

In this article, we explore the behavior of the sandpile models on evolutive graphs, that is, graphs that change their topology over time.
Evolutive sandpile had the potential to enrich the classic sandpile models.
It is worth mentioning for instance the extension of the Schelling segregation model over evolutive graphs in \cite{henry2011emergence}.
Also, Alan Kirman has proposed that economics could be effectively modeled as a network evolving over time, wherein agents have the ability to learn from past experiences with their neighbors, allowing links between them to either weaken or strengthen~\cite{kirman1997economy}. 
In this setting, it is possible to explore classic results and extend them to evolutive sandpiles.
We call them \emph{evolutive sandpiles} and define them in Section~\ref{sec:evosan}.
This turns out in the occurrence of phenomena impossible in the classical sandpile models.
For instance, configurations over evolutive graphs that are always unstable.
In Section~\ref{sec:patterns}, we explore the patterns obtained from stabilizing a starting configuration in which only one vertex contains grains.
In the classical sandpile models, interesting fractals emerge from these experiments.
In Section~\ref{sec:powerlaws}, we investigate the power laws that better fit the number of topplings required to stabilize evolutive sandpile models.

\section{Evolutive sandpiles}\label{sec:evosan}

Several different models of sandpiles have appeared, an account might be found at \cite{dhar1999abelian}.
A novel extension of the sandpile models occurs by changing the topology of the network over time, {\it i.e.}, the edges of the graph are allowed to change over discrete time, therefore the Laplacian matrix also will change over discrete time.
To make this more precise, let us introduce the {\it evolutive Laplacian matrix} $L(t)$ as follows.
For $t\in \mathbb{N}$, the $u,v$-entry of the evolutive Laplacian matrix is given by
\[
L(t)_{u,v}=
\begin{cases}
\sum_{w\in N(u)} f_{u,w}(t) & \text{if } u=v,\\
-f_{u,v}(t) & \text{otherwise},
\end{cases}
\]
where $f_{u,v}(t)$ is the function that provides the number of edges between vertices $u$ and $v$ at time $t$.
Therefore, for this evolutive Laplacian matrix, there is an associated an \textit{evolutive graph} $G_t$ at time $t$ which has edge set $E_t$.

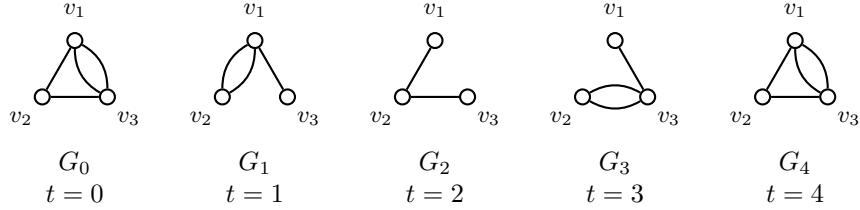
\begin{figure}[!ht]
    \centering
    \begin{tabular}{ccccc}
        \begin{tikzpicture}[scale=.5]
	        \tikzstyle{every node}=[minimum width=0pt, inner sep=2pt, circle]
            \draw (90:1) node (v1) [draw, label=above:{\small $v_1$}] {};
    		\draw (210:1) node (v2) [draw, label=below left:{\small $v_2$}] {};
    		\draw (330:1) node (v3) [draw, label=below right:{\small $v_3$}] {};
    		\draw (v1) -- (v2);
    		\draw (v1) edge[bend left] (v3);
    		\draw (v1) edge[bend right] (v3);
    		\draw (v3) -- (v2);
	    \end{tikzpicture}
        &
        \begin{tikzpicture}[scale=.5]
	        \tikzstyle{every node}=[minimum width=0pt, inner sep=2pt, circle]
            \draw (90:1) node (v1) [draw, label=above:{\small $v_1$}] {};
    		\draw (210:1) node (v2) [draw, label=below left:{\small $v_2$}] {};
    		\draw (330:1) node (v3) [draw, label=below right:{\small $v_3$}] {};
    		\draw (v1) edge[bend left] (v2);
    		\draw (v1) edge[bend right] (v2);
    		\draw (v1) -- (v3);
	    \end{tikzpicture}
        &
        \begin{tikzpicture}[scale=.5]
	        \tikzstyle{every node}=[minimum width=0pt, inner sep=2pt, circle]
            \draw (90:1) node (v1) [draw, label=above:{\small $v_1$}] {};
    		\draw (210:1) node (v2) [draw, label=below left:{\small $v_2$}] {};
    		\draw (330:1) node (v3) [draw, label=below right:{\small $v_3$}] {};
    		\draw (v2) -- (v1);
    		\draw (v2) -- (v3);
	    \end{tikzpicture}
        &
        \begin{tikzpicture}[scale=.5]
	        \tikzstyle{every node}=[minimum width=0pt, inner sep=2pt, circle]
            \draw (90:1) node (v1) [draw, label=above:{\small $v_1$}] {};
    		\draw (210:1) node (v2) [draw, label=below left:{\small $v_2$}] {};
    		\draw (330:1) node (v3) [draw, label=below right:{\small $v_3$}] {};
    		\draw (v1) -- (v3);
    		\draw (v2) edge[bend left] (v3);
    		\draw (v2) edge[bend right] (v3);
	    \end{tikzpicture}
	    &
	    \begin{tikzpicture}[scale=.5]
	        \tikzstyle{every node}=[minimum width=0pt, inner sep=2pt, circle]
            \draw (90:1) node (v1) [draw, label=above:{\small $v_1$}] {};
    		\draw (210:1) node (v2) [draw, label=below left:{\small $v_2$}] {};
    		\draw (330:1) node (v3) [draw, label=below right:{\small $v_3$}] {};
    		\draw (v1) -- (v2);
    		\draw (v1) edge[bend left] (v3);
    		\draw (v1) edge[bend right] (v3);
    		\draw (v3) -- (v2);
	    \end{tikzpicture}
        \\
        $G_0$ & $G_1$ & $G_2$ & $G_3$ & $G_4$
        \\
        $t=0$ & $t=1$ & $t=2$ & $t=3$ & $t=4$
    \end{tabular}
    \caption{Graph evolving over time $t=0,1,2,3,4$.}
    \label{fig:evolvinggraph}
\end{figure}

Consider the following evolutive Laplacian matrix


\[
\begin{bmatrix}
\sin\left(\frac{\pi t}{2}\right) + \cos\left(\frac{\pi t}{2}\right) + 2 & -\sin\left(\frac{\pi t}{2}\right) - 1 & -\cos\left(\frac{\pi t}{2}\right) - 1\\
 -\sin\left(\frac{\pi t}{2}\right) - 1 &  2 & \sin\left(\frac{\pi t}{2}\right) - 1\\
 -\cos\left(\frac{\pi t}{2}\right) - 1 & \sin\left(\frac{\pi t}{2}\right) - 1 & \cos\left(\frac{\pi t}{2}\right) - \sin\left(\frac{\pi t}{2}\right) + 2
\end{bmatrix}.
\]

Figure~\ref{fig:evolvinggraph} shows the evolution of the graph over time.

Then, the dynamics of this evolutive sandpile model start at time $t=0$, with a configuration vector $c$ that provides the number of grains of sand at each vertex of $G_0$.
A vertex $v\in G_t$ is {\it unstable} at time $t$ if $c_v\geq d_{G_t}(v)$.
Then, at time $t$ each unstable non-sink vertex is toppled once, and the obtained configuration will be the configuration at the following time $t+1$. 
Note that the configuration obtained after toppling each unstable vertex of $G_t$ is unique, this is due to the \emph{commutativity of the toppling process} see \cite[Theorem 2.2.2]{MR3889995}.

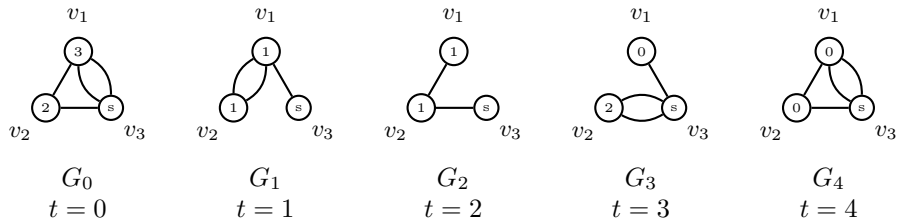
\begin{figure}[ht]
    \centering
    \begin{tabular}{ccccc}
        \begin{tikzpicture}[scale=.5]
	        \tikzstyle{every node}=[minimum width=0pt, inner sep=2pt, circle]
            \draw (90:1) node (v1) [draw, label=above:{\small $v_1$}] {\tiny 3};
    		\draw (210:1) node (v2) [draw, label=below left:{\small $v_2$}] {\tiny  2};
    		\draw (330:1) node (v3) [draw, label=below right:{\small $v_3$}] {\tiny  s};
    		\draw (v1) -- (v2);
    		\draw (v1) edge[bend left] (v3);
    		\draw (v1) edge[bend right] (v3);
    		\draw (v3) -- (v2);
	    \end{tikzpicture}
        &
        \begin{tikzpicture}[scale=.5]
	        \tikzstyle{every node}=[minimum width=0pt, inner sep=2pt, circle]
            \draw (90:1) node (v1) [draw, label=above:{\small $v_1$}] {\tiny 1};
    		\draw (210:1) node (v2) [draw, label=below left:{\small $v_2$}] {\tiny 1};
    		\draw (330:1) node (v3) [draw, label=below right:{\small $v_3$}] {\tiny s};
    		\draw (v1) edge[bend left] (v2);
    		\draw (v1) edge[bend right] (v2);
    		\draw (v1) -- (v3);
	    \end{tikzpicture}
        &
        \begin{tikzpicture}[scale=.5]
	        \tikzstyle{every node}=[minimum width=0pt, inner sep=2pt, circle]
            \draw (90:1) node (v1) [draw, label=above:{\small $v_1$}] {\tiny 1};
    		\draw (210:1) node (v2) [draw, label=below left:{\small $v_2$}] {\tiny 1};
    		\draw (330:1) node (v3) [draw, label=below right:{\small $v_3$}] {\tiny s};
    		\draw (v2) -- (v1);
    		\draw (v2) -- (v3);
	    \end{tikzpicture}
        &
        \begin{tikzpicture}[scale=.5]
	        \tikzstyle{every node}=[minimum width=0pt, inner sep=2pt, circle]
            \draw (90:1) node (v1) [draw, label=above:{\small $v_1$}] {\tiny 0};
    		\draw (210:1) node (v2) [draw, label=below left:{\small $v_2$}] {\tiny 2};
    		\draw (330:1) node (v3) [draw, label=below right:{\small $v_3$}] {\tiny s};
    		\draw (v1) -- (v3);
    		\draw (v2) edge[bend left] (v3);
    		\draw (v2) edge[bend right] (v3);
	    \end{tikzpicture}
	    &
	    \begin{tikzpicture}[scale=.5]
	        \tikzstyle{every node}=[minimum width=0pt, inner sep=2pt, circle]
            \draw (90:1) node (v1) [draw, label=above:{\small $v_1$}] {\tiny 0};
    		\draw (210:1) node (v2) [draw, label=below left:{\small $v_2$}] {\tiny 0};
    		\draw (330:1) node (v3) [draw, label=below right:{\small $v_3$}] {\tiny s};
    		\draw (v1) -- (v2);
    		\draw (v1) edge[bend left] (v3);
    		\draw (v1) edge[bend right] (v3);
    		\draw (v3) -- (v2);
	    \end{tikzpicture}
        \\
        $G_0$ & $G_1$ & $G_2$ & $G_3$ & $G_4$
        \\
        $t=0$ & $t=1$ & $t=2$ & $t=3$ & $t=4$
    \end{tabular}
    \caption{Graph evolving over time $t=0,1,2,3,4$.}
    \label{fig:evolvingconfiguration}
\end{figure}

Continuing the previous example, assume vertex $v_3$ is the sink vertex.
Consider the initial configuration $c=(3,2)$.
In Figure~\ref{fig:evolvingconfiguration}, we can observe an \textit{avalanche}, that is, a sequence of topplings of unstable non-sink vertices and evolutions of the graph.

\begin{figure}[ht]
    \centering
    \begin{tabular}{cc}
        \begin{tikzpicture}[scale=.5]
	        \tikzstyle{every node}=[minimum width=0pt, inner sep=2pt, circle]
            \draw (90:1) node (u) [draw, label=above:{\small $u$}] {\tiny 1};
    		\draw (210:1) node (v) [draw, label=below left:{\small $v$}] {\tiny  0};
    		\draw (330:1) node (s) [draw, label=below right:{\small $s$}] {\tiny  s};
    		\draw (u) -- (v);
    		\draw (v) -- (s);
	    \end{tikzpicture}
        &
        \begin{tikzpicture}[scale=.5]
	        \tikzstyle{every node}=[minimum width=0pt, inner sep=2pt, circle]
            \draw (90:1) node (u) [draw, label=above:{\small $u$}] {\tiny 0};
    		\draw (210:1) node (v) [draw, label=below left:{\small $v$}] {\tiny  1};
    		\draw (330:1) node (s) [draw, label=below right:{\small $s$}] {\tiny  s};
    		\draw (u) -- (v);
    		\draw (u) -- (s);
	    \end{tikzpicture}
        \\
        $G_{even}$ & $G_{odd}$
        \\
        $t=even$ & $t=odd$
    \end{tabular}
    \caption{An avalanche of configurations over an evolutive graph that never stabilizes.}
    \label{fig:evolvingconfigurationthatdoesnotstabilizes}
\end{figure}
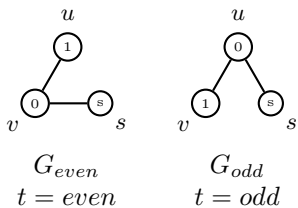

As previously mentioned, in the classic sandpile model, it is possible to show that for any finite graph $G$ with sink $s$, every configuration on $(G,s)$ eventually reaches a stable configuration under the toppling operations.
However, this is not always true over an evolutive graph with a sink.
For instance, consider an evolutive graph, whose vertex set is $\{u,v,s\}$ where $s$ is the sink, that evolves on two stages where $E_t=E(G_t)={uv,us}$ when $t$ is odd, and $E_t=E(G_t)={uv,vs}$ when $t$ is even.
Start at time $t=0$ with configuration $(1,0)$, where the first entry is associated with vertex $u$ and the second with the vertex $v$.
Since $\deg_{G_0}(u)=1$, then $u$ is unstable and sends its only grain of sand to $v$.
Then, at time $t=1$, the configuration is $(0,1)$, and vertex $v$ is unstable, since at this stage $\deg_{G_1}(v)=1$. 
Therefore, the configuration is $(1,0)$ when $t$ is even, and it is $(0,1)$ when $t$ is odd.
Figure~\ref{fig:evolvingconfigurationthatdoesnotstabilizes} shows the evolution of these configurations.
Is it possible to obtain a criterion to decide whether a configuration over an evolutive graph would stabilize after a finite number of topplings?

We continue the exploration of the notion of evolutive sandpiles under different perspectives in the following sections.

\section{Pattern formation}\label{sec:patterns}

\begin{figure}[ht!]
    \centering
    \begin{tabular}{c@{\extracolsep{1cm}}c}
    \begin{tabular}{ccc}
        \begin{tikzpicture}[scale=.5]
	        \tikzstyle{every node}=[minimum width=0pt, inner sep=2pt, circle]
            \foreach \i in {-1,...,1} {
                \draw [thick,draw] (\i,-1.3) -- (\i,1.3); 
            }
            \foreach \i in {-1,...,1} {
                \draw [thick,draw] (-1.3,\i) -- (1.3,\i); 
            }
            \draw (-1,1) node (v11) [draw,fill=white] {\tiny 0};
            \draw (-1,0) node (v12) [draw,fill=white] {\tiny 0};
            \draw (-1,-1) node (v13) [draw,fill=white] {\tiny 0};
            \draw (0,1) node (v21) [draw,fill=white] {\tiny 0};
            \draw (0,0) node (v22) [draw,fill=white] {\tiny 8};
            \draw (0,-1) node (v23) [draw,fill=white] {\tiny 0};
            \draw (1,1) node (v31) [draw,fill=white] {\tiny 0};
            \draw (1,0) node (v32) [draw,fill=white] {\tiny 0};
            \draw (1,-1) node (v33) [draw,fill=white] {\tiny 0};
            \draw (0,-2) node {$G_0$};
            \draw (0,-3) node {$t=0$};
	    \end{tikzpicture}
        &
        \begin{tikzpicture}[scale=.5]
	        \tikzstyle{every node}=[minimum width=0pt, inner sep=2pt, circle]
            \foreach \i in {-1,...,1} {
                \draw [thick,draw] (\i,-1.3) -- (\i,1.3); 
            }
            \foreach \i in {-1,...,1} {
                \draw [thick,draw] (-1.3,\i) -- (1.3,\i); 
            }
            \draw (-1,1) node (v11) [draw,fill=white] {\tiny 0};
            \draw (-1,0) node (v12) [draw,fill=white] {\tiny 1};
            \draw (-1,-1) node (v13) [draw,fill=white] {\tiny 0};
            \draw (0,1) node (v21) [draw,fill=white] {\tiny 1};
            \draw (0,0) node (v22) [draw,fill=white] {\tiny 4};
            \draw (0,-1) node (v23) [draw,fill=white] {\tiny 1};
            \draw (1,1) node (v31) [draw,fill=white] {\tiny 0};
            \draw (1,0) node (v32) [draw,fill=white] {\tiny 1};
            \draw (1,-1) node (v33) [draw,fill=white] {\tiny 0};
            \draw (0,-2) node {$G_1$};
            \draw (0,-3) node {$t=1$};
	    \end{tikzpicture}
        &
        \begin{tikzpicture}[scale=.5]
	        \tikzstyle{every node}=[minimum width=0pt, inner sep=2pt, circle]
            \foreach \i in {-1,...,1} {
                \draw [thick,draw] (\i,-1.3) -- (\i,1.3); 
            }
            \foreach \i in {-1,...,1} {
                \draw [thick,draw] (-1.3,\i) -- (1.3,\i); 
            }
            \draw (-1,1) node (v11) [draw,fill=white] {\tiny 0};
            \draw (-1,0) node (v12) [draw,fill=white] {\tiny 2};
            \draw (-1,-1) node (v13) [draw,fill=white] {\tiny 0};
            \draw (0,1) node (v21) [draw,fill=white] {\tiny 2};
            \draw (0,0) node (v22) [draw,fill=white] {\tiny 0};
            \draw (0,-1) node (v23) [draw,fill=white] {\tiny 2};
            \draw (1,1) node (v31) [draw,fill=white] {\tiny 0};
            \draw (1,0) node (v32) [draw,fill=white] {\tiny 2};
            \draw (1,-1) node (v33) [draw,fill=white] {\tiny 0};
            \draw (0,-2) node {$G_2$};
            \draw (0,-3) node {$t=2$};
	    \end{tikzpicture}
         
    \end{tabular}

    &
    \raisebox{-.5\height}{\includegraphics[width=4cm]{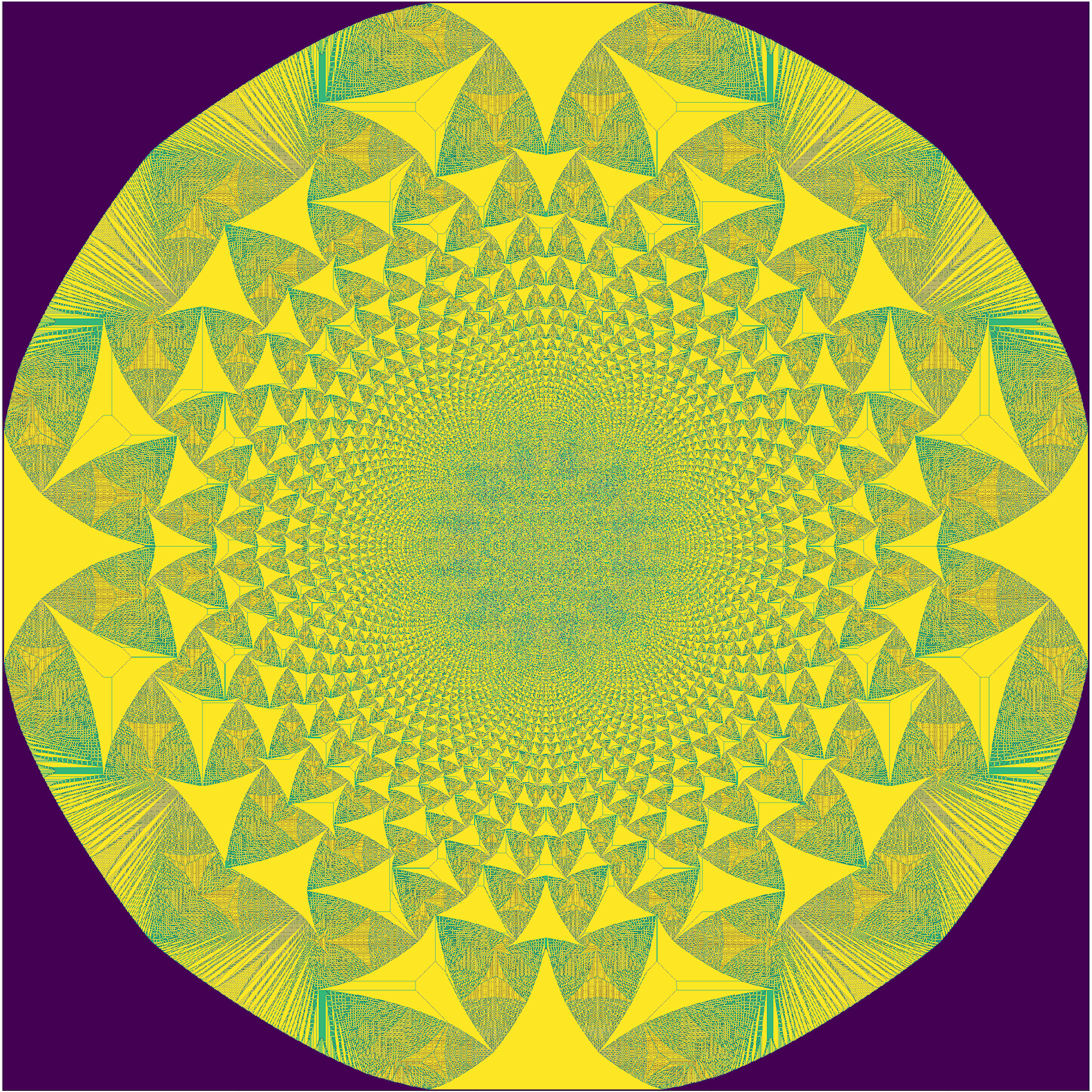}}
        \\
        (a) & (b)\\
    \begin{tabular}{ccc}
        \begin{tikzpicture}[scale=.5]
	        \tikzstyle{every node}=[minimum width=0pt, inner sep=2pt, circle]
            \foreach \i in {-1,...,1} {
                \draw [thick,draw] (\i,-1.3) -- (\i,1.3); 
            }
            \foreach \i in {-1,...,1} {
                \draw [thick,draw] (-1.3,\i) -- (1.3,\i); 
            }
            \draw (-1,1) node (v11) [draw,fill=white] {\tiny 0};
            \draw (-1,0) node (v12) [draw,fill=white] {\tiny 0};
            \draw (-1,-1) node (v13) [draw,fill=white] {\tiny 0};
            \draw (0,1) node (v21) [draw,fill=white] {\tiny 0};
            \draw (0,0) node (v22) [draw,fill=white] {\tiny 8};
            \draw (0,-1) node (v23) [draw,fill=white] {\tiny 0};
            \draw (1,1) node (v31) [draw,fill=white] {\tiny 0};
            \draw (1,0) node (v32) [draw,fill=white] {\tiny 0};
            \draw (1,-1) node (v33) [draw,fill=white] {\tiny 0};
            \draw (0,-2) node {$G_0$};
            \draw (0,-3) node {$t=0$};
	    \end{tikzpicture}
        &
        \begin{tikzpicture}[scale=.5]
	        \tikzstyle{every node}=[minimum width=0pt, inner sep=2pt, circle]
            \draw (-1,1) node (v11) [draw,fill=white] {\tiny 0};
            \draw (-1,0) node (v12) [draw,fill=white] {\tiny 1};
            \draw (-1,-1) node (v13) [draw,fill=white] {\tiny 0};
            \draw (0,1) node (v21) [draw,fill=white] {\tiny 1};
            \draw (0,0) node (v22) [draw,fill=white] {\tiny 4};
            \draw (0,-1) node (v23) [draw,fill=white] {\tiny 1};
            \draw (1,1) node (v31) [draw,fill=white] {\tiny 0};
            \draw (1,0) node (v32) [draw,fill=white] {\tiny 1};
            \draw (1,-1) node (v33) [draw,fill=white] {\tiny 0};
            \draw (0,-2) node {$G_1$};
            \draw (0,-3) node {$t=1$};
            \draw (v11) -- (v22) -- (v33);
            \draw (v13) -- (v22) -- (v31);
            \draw (v12) -- (v21) -- (v32) -- (v23) -- (v12);
	    \end{tikzpicture}
        &
        \begin{tikzpicture}[scale=.5]
	        \tikzstyle{every node}=[minimum width=0pt, inner sep=2pt, circle]
            \foreach \i in {-1,...,1} {
                \draw [thick,draw] (\i,-1.3) -- (\i,1.3); 
            }
            \foreach \i in {-1,...,1} {
                \draw [thick,draw] (-1.3,\i) -- (1.3,\i); 
            }
            \draw (-1,1) node (v11) [draw,fill=white] {\tiny 1};
            \draw (-1,0) node (v12) [draw,fill=white] {\tiny 1};
            \draw (-1,-1) node (v13) [draw,fill=white] {\tiny 1};
            \draw (0,1) node (v21) [draw,fill=white] {\tiny 1};
            \draw (0,0) node (v22) [draw,fill=white] {\tiny 0};
            \draw (0,-1) node (v23) [draw,fill=white] {\tiny 1};
            \draw (1,1) node (v31) [draw,fill=white] {\tiny 1};
            \draw (1,0) node (v32) [draw,fill=white] {\tiny 1};
            \draw (1,-1) node (v33) [draw,fill=white] {\tiny 1};
            \draw (0,-2) node {$G_2$};
            \draw (0,-3) node {$t=2$};
	    \end{tikzpicture}
    \end{tabular}
    
        &
        \raisebox{-.5\height}{\includegraphics[width=4cm]{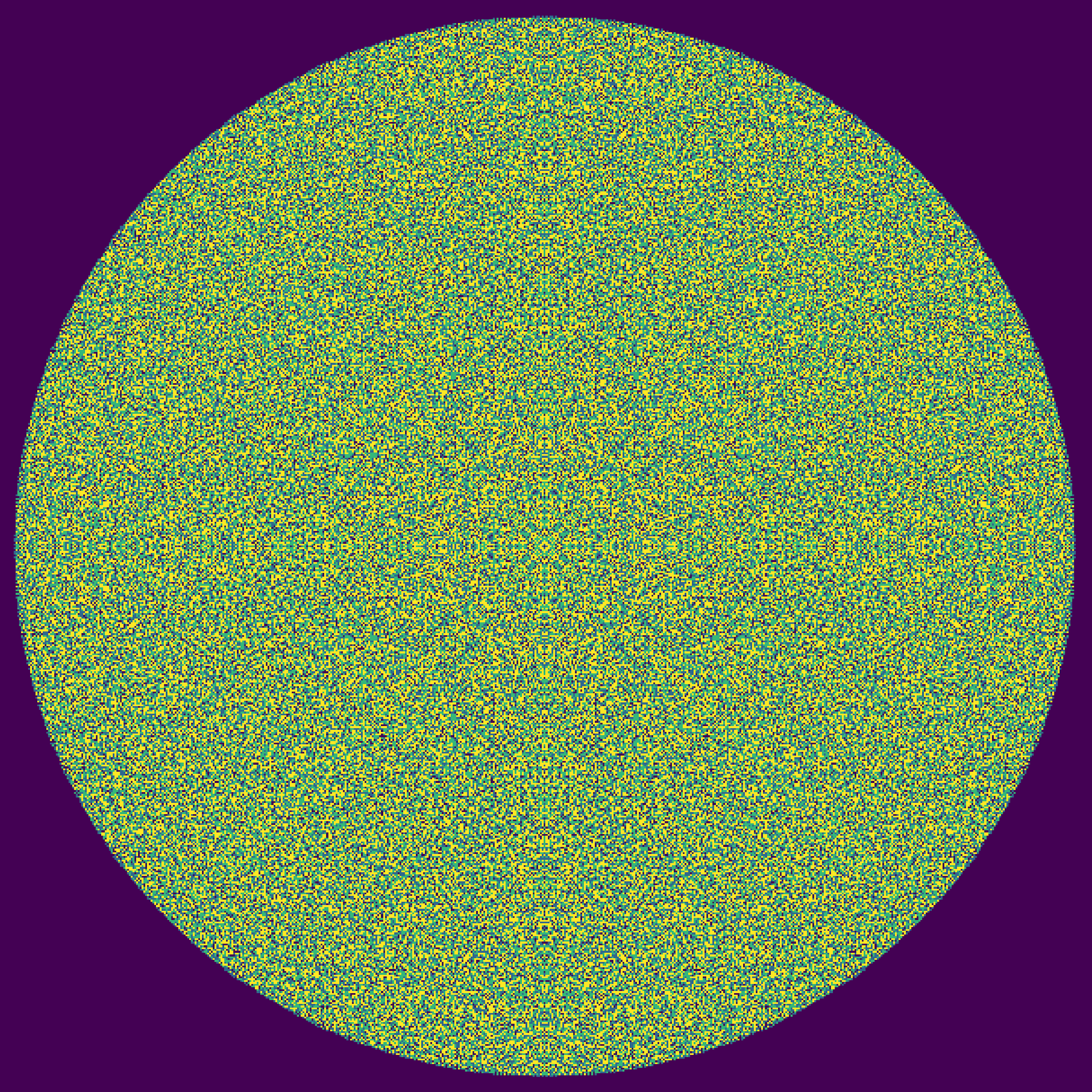}}
        \\
        (c) & (d) \\
        \begin{tabular}{cc}
             \begin{tikzpicture}[scale=.5]
	        \tikzstyle{every node}=[minimum width=0pt, inner sep=2pt, circle]
            \foreach \i in {-1,...,1} {
                \foreach \j in {-1,0}{
                    \draw [thick,draw] (\i,\j) edge[bend right] (\i,\j+1); 
                    \draw [thick,draw] (\i,\j) edge[bend left] (\i,\j+1);
                }
            }
            \foreach \i in {-1,...,1} {
                \draw [thick,draw] (-1.3,\i) -- (1.3,\i); 
            }
            \draw (-1,1) node (v11) [draw,fill=white] {\tiny 0};
            \draw (-1,0) node (v12) [draw,fill=white] {\tiny 0};
            \draw (-1,-1) node (v13) [draw,fill=white] {\tiny 0};
            \draw (0,1) node (v21) [draw,fill=white] {\tiny 0};
            \draw (0,0) node (v22) [draw,fill=white] {\tiny 8};
            \draw (0,-1) node (v23) [draw,fill=white] {\tiny 0};
            \draw (1,1) node (v31) [draw,fill=white] {\tiny 0};
            \draw (1,0) node (v32) [draw,fill=white] {\tiny 0};
            \draw (1,-1) node (v33) [draw,fill=white] {\tiny 0};
            \draw (0,-2) node {$G_0$};
            \draw (0,-3) node {$t=0$};
	    \end{tikzpicture}
        &
        \begin{tikzpicture}[scale=.5]
	        \tikzstyle{every node}=[minimum width=0pt, inner sep=2pt, circle]
            \foreach \i in {-1,...,1} {
                \foreach \j in {-1,0}{
                    \draw [thick,draw] (\j,\i) edge[bend right] (\j+1,\i); 
                    \draw [thick,draw] (\j,\i) edge[bend left] (\j+1,\i);
                }
            }
            \foreach \i in {-1,...,1} {
                \draw [thick,draw] (\i,-1.3) -- (\i,1.3); 
            }
            \draw (-1,1) node (v11) [draw,fill=white] {\tiny 0};
            \draw (-1,0) node (v12) [draw,fill=white] {\tiny 1};
            \draw (-1,-1) node (v13) [draw,fill=white] {\tiny 0};
            \draw (0,1) node (v21) [draw,fill=white] {\tiny 2};
            \draw (0,0) node (v22) [draw,fill=white] {\tiny 2};
            \draw (0,-1) node (v23) [draw,fill=white] {\tiny 2};
            \draw (1,1) node (v31) [draw,fill=white] {\tiny 0};
            \draw (1,0) node (v32) [draw,fill=white] {\tiny 1};
            \draw (1,-1) node (v33) [draw,fill=white] {\tiny 0};
            \draw (0,-2) node {$G_1$};
            \draw (0,-3) node {$t=1$};
	    \end{tikzpicture}
        \end{tabular}
        &
        \raisebox{-.5\height}{\includegraphics[width=4cm]{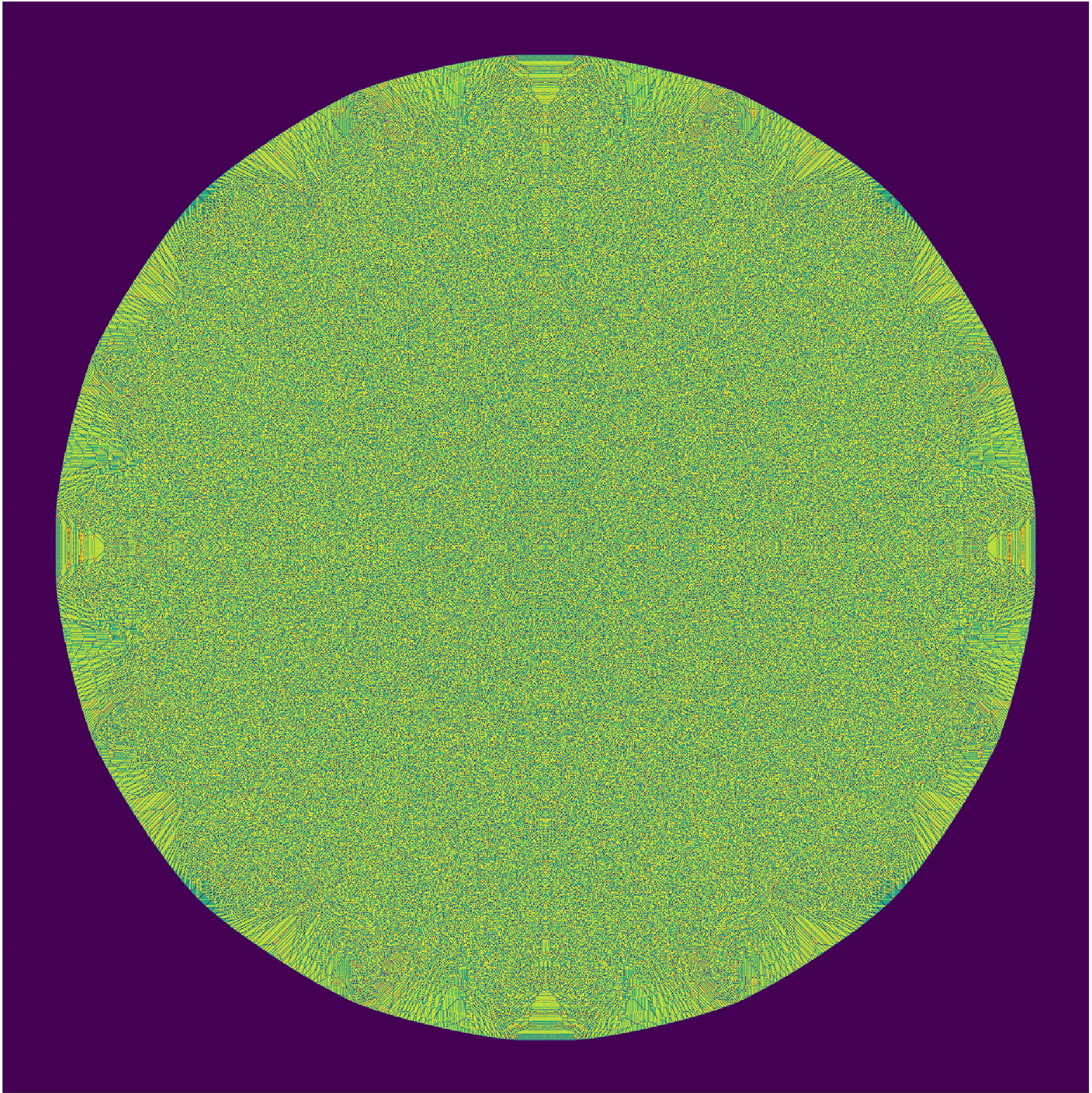}}
        \\
        (e) & (f)
    \end{tabular}
    
    \caption{(a) (c) and (e) are partial views of infinite evolutive graphs of stabilizations of an initial configuration having $8$ grains at the center evolving with time $t$. 
    In (a) the evolutive graph is the same at any time. 
    (b) The stabilization of the configuration with $10,000,000$ grains at the center over the same evolutive graph as in (a). 
    In (c) and (e), the evolutive graphs evolve only over two different infinite graphs. 
    (d) The stabilization of the configuration with $500,000$ grains at the center over the same evolutive graph as in (c).
    (f) The stabilization of the configuration with $5,000,000$ grains at the center over the same evolutive graph as in (e).
    }
    \label{fig:evolvingsandpilegrid}
\end{figure}
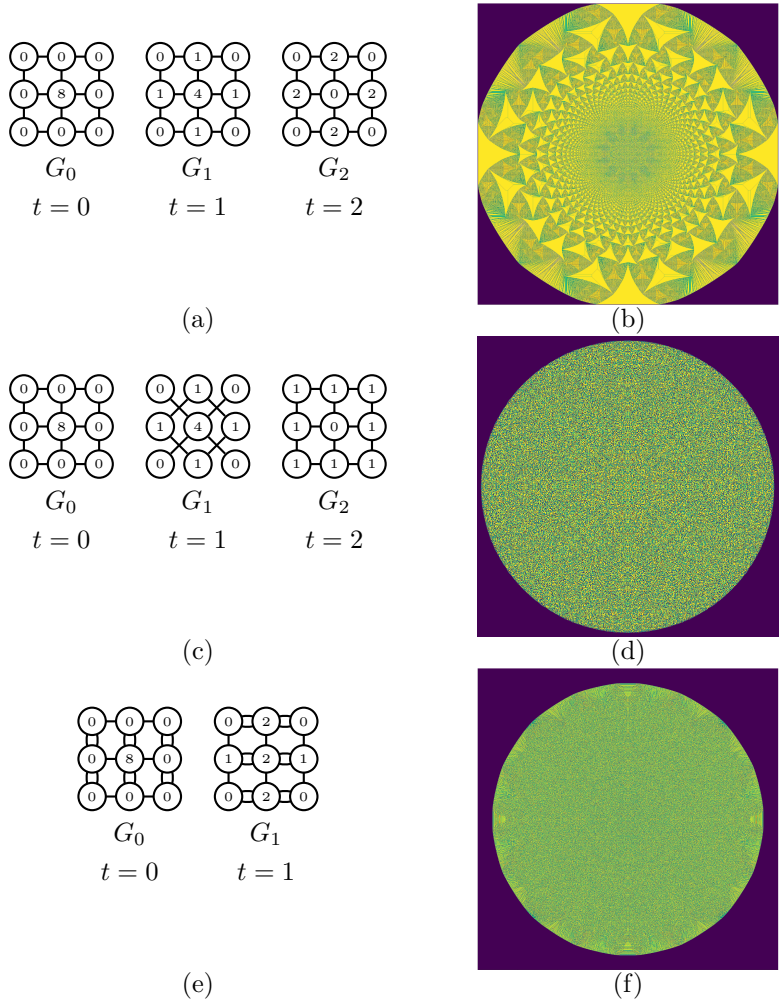

Many geometric structures have appeared within the context of sandpiles, most of them exhibiting interesting fractal patterns. 
One such geometric structure is obtained on the infinite square grid graph, where the starting configuration has a large number of grains at the center vertex $(0,0)$ and no grains elsewhere. 
After stabilizing this configuration, a fractal starts to emerge, see Fig.~\ref{fig:evolvingsandpilegrid} (b).
It turns out that the continuum scaling limit of this Abelian sandpile has the structure of an Apollonian circle packing, the details are out of the scope of this text for which the reader, interested in deeply understanding the mathematics behind this fractal, is invited to see~\cite{MR3494492,MR3664999}.
Many variants have appeared since then, for example, several authors have studied the patterns obtained when the initial configuration has a fixed height at each non-center vertex. 
The geometry when the initial height is equal to 2 at the non-center vertices is investigated in \cite{OSTOJIC2003187}.

We explore the behavior of the stabilization of a configuration over an evolutive graph in which a unique vertex is featured with a large number of grains of sand.
In these experiments, the evolutive graph has an infinite number of vertices but the number of edges joining any pair of vertices is finite.
First, we observe that classic fractals can be obtained by using the same graph over time, see Figure~\ref{fig:evolvingsandpilegrid} (a) and (b).
In this figure, at (a) we see the stabilization of the all-zeros configuration except by 8 grains at the center of the infinite grid that remains the same over time. Note that the stabilization occurs only on the 3-by-3 sub-grid for that reason we only show this sub-grid.
At (b), we stabilize the starting configuration with the all-zeros at each vertex except by 10,000,000 grains placed at the center of the infinite grid.
At (c) and (d), the evolutive graph only has two stages, even cases correspond to the infinite square grid, and the odd cases correspond to the infinite grid in which each vertex at $(i,j)$ is adjacent only with vertices at sites $(i-1,j+1)$, $(i+1,j+1)$, $(i-1,j-1)$ and $(i+1,j-1)$.
Note the graph in odd cases has two connected components.
At (c), we see the stabilization of the all-zeros configuration except by 8 grains at the center over the described evolutive graph.
And at (d), we have the stabilization of the configuration with the all-zeros at each vertex except by 500,000 grains placed at the center.
It is interesting to observe that no fractal emerges in this evolutive sandpile. However, some symmetry over the horizontal and vertical axes might be identified. 
At (e) and (f), the evolutive graph only has two stages, even cases correspond to the infinite square grid adding additional vertical edges, and the odd cases correspond to the infinite square grid adding additional horizontal edges.
In both cases, each note has degree 5.
At (e), we see the stabilization of the all-zeros configuration except by 8 grains at the center over the described evolutive graph.
Lastly, at (f), we have the stabilization of the configuration with the all-zeros at each vertex except by 5,000,000 grains placed at the center.

\begin{figure}[ht!]
    \centering
    \begin{tabular}{cc}

        \begin{tabular}{ccc}
        \begin{tikzpicture}[scale=.5]
	        \tikzstyle{every node}=[minimum width=0pt, inner sep=2pt, circle]
            \foreach \i in {-1,...,1} {
                \draw [thick,draw] (\i,-1.3) -- (\i,1.3); 
            }
            \draw (-1,1) node (v11) [draw,fill=white] {\tiny 0};
            \draw (-1,0) node (v12) [draw,fill=white] {\tiny 0};
            \draw (-1,-1) node (v13) [draw,fill=white] {\tiny 0};
            \draw (0,1) node (v21) [draw,fill=white] {\tiny 0};
            \draw (0,0) node (v22) [draw,fill=white] {\tiny 8};
            \draw (0,-1) node (v23) [draw,fill=white] {\tiny 0};
            \draw (1,1) node (v31) [draw,fill=white] {\tiny 0};
            \draw (1,0) node (v32) [draw,fill=white] {\tiny 0};
            \draw (1,-1) node (v33) [draw,fill=white] {\tiny 0};
	    \end{tikzpicture}
        &
        \begin{tikzpicture}[scale=.5]
	        \tikzstyle{every node}=[minimum width=0pt, inner sep=2pt, circle]
            
            \foreach \i in {-1,...,1} {
                \draw [thick,draw] (-1.3,\i) -- (1.3,\i); 
            }
            \draw (-1,1) node (v11) [draw,fill=white] {\tiny 0};
            \draw (-1,0) node (v12) [draw,fill=white] {\tiny 0};
            \draw (-1,-1) node (v13) [draw,fill=white] {\tiny 0};
            \draw (0,1) node (v21) [draw,fill=white] {\tiny 1};
            \draw (0,0) node (v22) [draw,fill=white] {\tiny 6};
            \draw (0,-1) node (v23) [draw,fill=white] {\tiny 1};
            \draw (1,1) node (v31) [draw,fill=white] {\tiny 0};
            \draw (1,0) node (v32) [draw,fill=white] {\tiny 0};
            \draw (1,-1) node (v33) [draw,fill=white] {\tiny 0};
	    \end{tikzpicture}
        &
        \begin{tikzpicture}[scale=.5]
	        \tikzstyle{every node}=[minimum width=0pt, inner sep=2pt, circle]
            \foreach \i in {-1,...,1} {
                \draw [thick,draw] (\i,-1.3) -- (\i,1.3); 
            }
            \draw (-1,1) node (v11) [draw,fill=white] {\tiny 0};
            \draw (-1,0) node (v12) [draw,fill=white] {\tiny 1};
            \draw (-1,-1) node (v13) [draw,fill=white] {\tiny 0};
            \draw (0,1) node (v21) [draw,fill=white] {\tiny 1};
            \draw (0,0) node (v22) [draw,fill=white] {\tiny 4};
            \draw (0,-1) node (v23) [draw,fill=white] {\tiny 1};
            \draw (1,1) node (v31) [draw,fill=white] {\tiny 0};
            \draw (1,0) node (v32) [draw,fill=white] {\tiny 1};
            \draw (1,-1) node (v33) [draw,fill=white] {\tiny 0};
	    \end{tikzpicture}
        \\
         $G_0$ & $G_1$ & $G_2$
        \\
        $t=0$ & $t=1$ & $t=2$
        \\
        \begin{tikzpicture}[scale=.5]
	        \tikzstyle{every node}=[minimum width=0pt, inner sep=2pt, circle]
            
            \foreach \i in {-1,...,1} {
                \draw [thick,draw] (-1.3,\i) -- (1.3,\i); 
            }
            \draw (-1,1) node (v11) [draw,fill=white] {\tiny 0};
            \draw (-1,0) node (v12) [draw,fill=white] {\tiny 1};
            \draw (-1,-1) node (v13) [draw,fill=white] {\tiny 0};
            \draw (0,1) node (v21) [draw,fill=white] {\tiny 2};
            \draw (0,0) node (v22) [draw,fill=white] {\tiny 2};
            \draw (0,-1) node (v23) [draw,fill=white] {\tiny 2};
            \draw (1,1) node (v31) [draw,fill=white] {\tiny 0};
            \draw (1,0) node (v32) [draw,fill=white] {\tiny 1};
            \draw (1,-1) node (v33) [draw,fill=white] {\tiny 0};
	    \end{tikzpicture}
	    &
	    \begin{tikzpicture}[scale=.5]
	        \tikzstyle{every node}=[minimum width=0pt, inner sep=2pt, circle]
            \foreach \i in {-1,...,1} {
                \draw [thick,draw] (\i,-1.3) -- (\i,1.3); 
            }
            \draw (-1,1) node (v11) [draw,fill=white] {\tiny 1};
            \draw (-1,0) node (v12) [draw,fill=white] {\tiny 2};
            \draw (-1,-1) node (v13) [draw,fill=white] {\tiny 1};
            \draw (0,1) node (v21) [draw,fill=white] {\tiny 0};
            \draw (0,0) node (v22) [draw,fill=white] {\tiny 0};
            \draw (0,-1) node (v23) [draw,fill=white] {\tiny 0};
            \draw (1,1) node (v31) [draw,fill=white] {\tiny 1};
            \draw (1,0) node (v32) [draw,fill=white] {\tiny 2};
            \draw (1,-1) node (v33) [draw,fill=white] {\tiny 1};
	    \end{tikzpicture}
	    &
	    \begin{tikzpicture}[scale=.5]
	        \tikzstyle{every node}=[minimum width=0pt, inner sep=2pt, circle]
            \foreach \i in {-1,...,1} {
                \draw [thick,draw] (-1.3,\i) -- (1.3,\i); 
            }
            \draw (-1,1) node (v11) [draw,fill=white] {\tiny 2};
            \draw (-1,0) node (v12) [draw,fill=white] {\tiny 0};
            \draw (-1,-1) node (v13) [draw,fill=white] {\tiny 2};
            \draw (0,1) node (v21) [draw,fill=white] {\tiny 0};
            \draw (0,0) node (v22) [draw,fill=white] {\tiny 0};
            \draw (0,-1) node (v23) [draw,fill=white] {\tiny 0};
            \draw (1,1) node (v31) [draw,fill=white] {\tiny 2};
            \draw (1,0) node (v32) [draw,fill=white] {\tiny 0};
            \draw (1,-1) node (v33) [draw,fill=white] {\tiny 2};
	    \end{tikzpicture}
        \\
        $G_3$ & $G_4$ & $G_5$
        \\
        $t=3$ & $t=4$ & $t=5$
        \\
        \begin{tikzpicture}[scale=.5]
	        \tikzstyle{every node}=[minimum width=0pt, inner sep=2pt, circle]
            \foreach \i in {-2,...,2} {
                \draw [thick,draw] (\i,-2.3) -- (\i,2.3); 
            }
            \draw (-2,2) node (v00) [draw,fill=white] {\tiny 0};
            \draw (-2,1) node (v01) [draw,fill=white] {\tiny 1};
            \draw (-2,0) node (v02) [draw,fill=white] {\tiny 0};
            \draw (-2,-1) node (v03) [draw,fill=white] {\tiny 1};
            \draw (-2,-2) node (v04) [draw,fill=white] {\tiny 0};
            \draw (-1,2) node (v10) [draw,fill=white] {\tiny 0};
            \draw (-1,1) node (v11) [draw,fill=white] {\tiny 0};
            \draw (-1,0) node (v12) [draw,fill=white] {\tiny 0};
            \draw (-1,-1) node (v13) [draw,fill=white] {\tiny 0};
            \draw (-1,-2) node (v14) [draw,fill=white] {\tiny 0};
            \draw (0,2) node (v20) [draw,fill=white] {\tiny 0};
            \draw (0,1) node (v21) [draw,fill=white] {\tiny 2};
            \draw (0,0) node (v22) [draw,fill=white] {\tiny 0};
            \draw (0,-1) node (v23) [draw,fill=white] {\tiny 2};
            \draw (0,-2) node (v24) [draw,fill=white] {\tiny 0};
            \draw (1,2) node (v30) [draw,fill=white] {\tiny 0};
            \draw (1,1) node (v31) [draw,fill=white] {\tiny 0};
            \draw (1,0) node (v32) [draw,fill=white] {\tiny 0};
            \draw (1,-1) node (v33) [draw,fill=white] {\tiny 0};
            \draw (1,-2) node (v34) [draw,fill=white] {\tiny 0};
            \draw (2,2) node (v40) [draw,fill=white] {\tiny 0};
            \draw (2,1) node (v41) [draw,fill=white] {\tiny 1};
            \draw (2,0) node (v42) [draw,fill=white] {\tiny 0};
            \draw (2,-1) node (v43) [draw,fill=white] {\tiny 1};
            \draw (2,-2) node (v44) [draw,fill=white] {\tiny 0};
	    \end{tikzpicture}
        &
        \begin{tikzpicture}[scale=.5]
	        \tikzstyle{every node}=[minimum width=0pt, inner sep=2pt, circle]
            \foreach \i in {-2,...,2} {
                \draw [thick,draw] (-2.3,\i) -- (2.3,\i); 
            }
            \draw (-2,2) node (v00) [draw,fill=white] {\tiny 0};
            \draw (-2,1) node (v01) [draw,fill=white] {\tiny 1};
            \draw (-2,0) node (v02) [draw,fill=white] {\tiny 0};
            \draw (-2,-1) node (v03) [draw,fill=white] {\tiny 1};
            \draw (-2,-2) node (v04) [draw,fill=white] {\tiny 0};
            \draw (-1,2) node (v10) [draw,fill=white] {\tiny 0};
            \draw (-1,1) node (v11) [draw,fill=white] {\tiny 0};
            \draw (-1,0) node (v12) [draw,fill=white] {\tiny 0};
            \draw (-1,-1) node (v13) [draw,fill=white] {\tiny 0};
            \draw (-1,-2) node (v14) [draw,fill=white] {\tiny 0};
            \draw (0,2) node (v20) [draw,fill=white] {\tiny 1};
            \draw (0,1) node (v21) [draw,fill=white] {\tiny 0};
            \draw (0,0) node (v22) [draw,fill=white] {\tiny 2};
            \draw (0,-1) node (v23) [draw,fill=white] {\tiny 0};
            \draw (0,-2) node (v24) [draw,fill=white] {\tiny 1};
            \draw (1,2) node (v30) [draw,fill=white] {\tiny 0};
            \draw (1,1) node (v31) [draw,fill=white] {\tiny 0};
            \draw (1,0) node (v32) [draw,fill=white] {\tiny 0};
            \draw (1,-1) node (v33) [draw,fill=white] {\tiny 0};
            \draw (1,-2) node (v34) [draw,fill=white] {\tiny 0};
            \draw (2,2) node (v40) [draw,fill=white] {\tiny 0};
            \draw (2,1) node (v41) [draw,fill=white] {\tiny 1};
            \draw (2,0) node (v42) [draw,fill=white] {\tiny 0};
            \draw (2,-1) node (v43) [draw,fill=white] {\tiny 1};
            \draw (2,-2) node (v44) [draw,fill=white] {\tiny 0};
	    \end{tikzpicture}
        &
        \begin{tikzpicture}[scale=.5]
	        \tikzstyle{every node}=[minimum width=0pt, inner sep=2pt, circle]
            \foreach \i in {-2,...,2} {
                \draw [thick,draw] (-2.3,\i) -- (2.3,\i); 
            }
            \draw (-2,2) node (v00) [draw,fill=white] {\tiny 0};
            \draw (-2,1) node (v01) [draw,fill=white] {\tiny 1};
            \draw (-2,0) node (v02) [draw,fill=white] {\tiny 0};
            \draw (-2,-1) node (v03) [draw,fill=white] {\tiny 1};
            \draw (-2,-2) node (v04) [draw,fill=white] {\tiny 0};
            \draw (-1,2) node (v10) [draw,fill=white] {\tiny 0};
            \draw (-1,1) node (v11) [draw,fill=white] {\tiny 0};
            \draw (-1,0) node (v12) [draw,fill=white] {\tiny 1};
            \draw (-1,-1) node (v13) [draw,fill=white] {\tiny 0};
            \draw (-1,-2) node (v14) [draw,fill=white] {\tiny 0};
            \draw (0,2) node (v20) [draw,fill=white] {\tiny 1};
            \draw (0,1) node (v21) [draw,fill=white] {\tiny 0};
            \draw (0,0) node (v22) [draw,fill=white] {\tiny 0};
            \draw (0,-1) node (v23) [draw,fill=white] {\tiny 0};
            \draw (0,-2) node (v24) [draw,fill=white] {\tiny 1};
            \draw (1,2) node (v30) [draw,fill=white] {\tiny 0};
            \draw (1,1) node (v31) [draw,fill=white] {\tiny 0};
            \draw (1,0) node (v32) [draw,fill=white] {\tiny 1};
            \draw (1,-1) node (v33) [draw,fill=white] {\tiny 0};
            \draw (1,-2) node (v34) [draw,fill=white] {\tiny 0};
            \draw (2,2) node (v40) [draw,fill=white] {\tiny 0};
            \draw (2,1) node (v41) [draw,fill=white] {\tiny 1};
            \draw (2,0) node (v42) [draw,fill=white] {\tiny 0};
            \draw (2,-1) node (v43) [draw,fill=white] {\tiny 1};
            \draw (2,-2) node (v44) [draw,fill=white] {\tiny 0};
	    \end{tikzpicture}
        \\
        $G_6$ & $G_7$ & $G_8$
        \\
        $t=6$ & $t=7$ & $t=8$
    \end{tabular}
    
         & \raisebox{-.5\height}{\includegraphics[width=4cm]{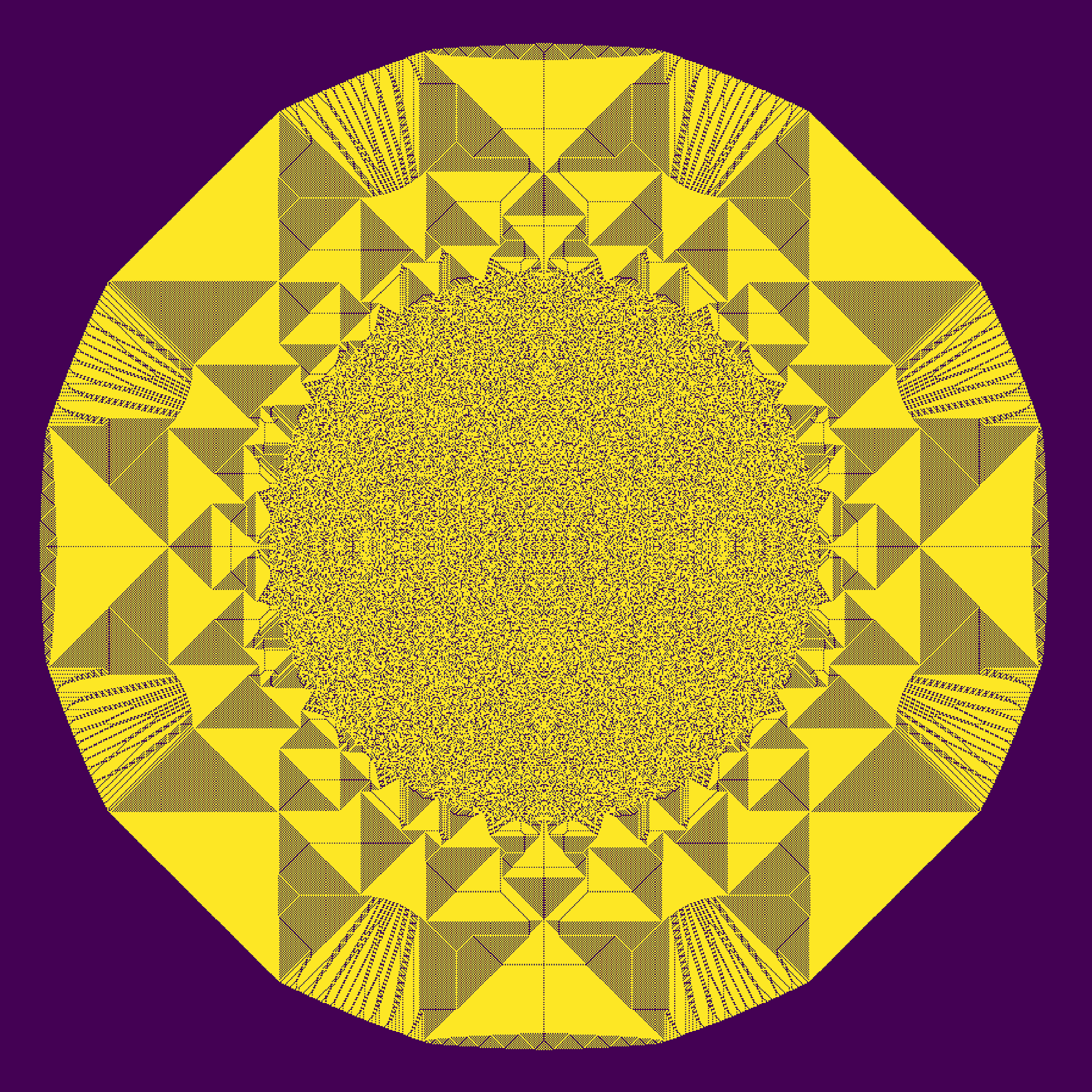}} \\
        (a) & (b)
    \end{tabular}

    \caption{(a) A partial view of the infinite evolutive graph of an stabilization of an evolutive sandpile where the first configuration has $8$ grains at the center evolving with time $t=0,1,2,3,4,5,6,7,8$. (b) The stabilization of the configuration with $400,000$ grains at the center over the same evolutive graph as in (a).}
    \label{fig:evolvingsandpilegrid3}
\end{figure}

In Figure~\ref{fig:evolvingsandpilegrid3}, we explore the evolutive graph consisting of two stages: the even case consists of the infinite graph in which vertex at site $(i,j)$ is adjacent only with vertices at sites $(i,j+1)$ and $(i,j-1)$, and the odd case consists of the infinite graph in which vertex at site $(i,j)$ is adjacent only with vertices at sites $(i+1,j)$ and $(i-1,j)$.
The topological properties of this evolutive graph consist in that each vertex has degree two and there is an infinite number of connected components at each stage.
At (a), we stabilize the configuration with only 8 grains at the center, this stabilization takes a total of 9 stages to get stabilized.
At (b), we have the stabilization of the configuration with $400,000$ grains at the center.
Changing the topology of the graph over time might result in new fractals and observations.

\section{Power laws}\label{sec:powerlaws}

Per Bak was a pioneer in the study of self-organized criticality. He observed that many natural processes occur not due to external forces, but through the evolution of the processes themselves. In these systems, events of any size can occur over time, and the process will self-organize until it reaches a stable state~\cite{bak1991catastrophes}. This is the central hypothesis of self-organized critical systems. Empirical evidence shows that earthquakes in a certain region of the Earth obey the Gutenberg-Richter law: the logarithm of the energy released by an earthquake is linearly related to the logarithm of the frequency of earthquakes of that energy. An earthquake continues until all the accumulated energy is released, reaching an energy-stable state. Therefore, earthquakes are examples of self-organizing systems. It is worth mentioning that in a self-organized system, the sizes of large and small events, are closely related by a scaling law. The number $N(x)$ of events of size $x$ is proportional to $x^{-\alpha}$ for some positive value $\alpha$. It is important to note that as the sizes of these events decrease, the number of events increases. This is because $N(x)$ tends to infinity as $x$ tends to zero. Therefore, it would not be appropriate to consider events of arbitrarily small size, but rather sizes above some cut-off value $x_{\text{min}}$. The archetypal model of self-organized criticality is the \emph{sandpile model}, in which a single grain of sand can trigger a large avalanche. The statistics behind such processes follow a power law, indicating a relationship between the number of large events and the number of small events. 

Now, our focus lies on sandpile models in which adjacencies evolve over time, investigating whether they adhere to a self-organized critical system or not. Similar to that of stationary systems. 
To accomplish this, we propose the following dynamics. Begin the evolutive sandpile model by randomly assigning a non-negative number, less than or equal to $\Delta(G_{0})-1$ to each vertex, where $G_{0}$ is the graph at time zero. This generates an initial random configuration. At time $t$, one grain is added to a randomly selected vertex, leading to the stabilization of $G_{t}$. 
This stabilization results in a new configuration $c'$, which serves as the input configuration for $G_{t+1}$, and this process repeats until the number of topplings becomes zero.

\begin{figure}
    \centering
    \begin{subfigure}{0.4\textwidth}
        \centering
        \includegraphics[width=\textwidth]{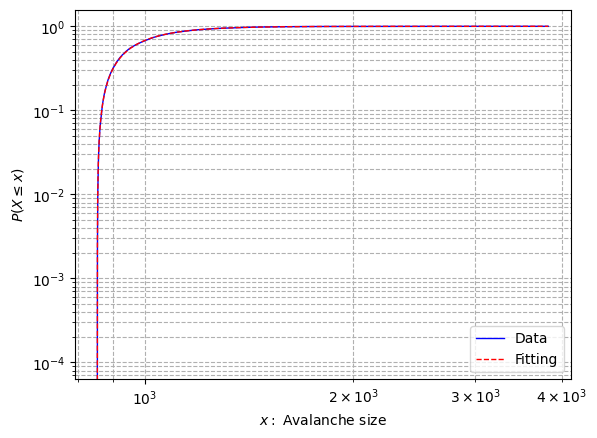}
        \caption{Model $\mathcal{G}$ of size $50\times 50$.}
        \label{fig:G50x50}
    \end{subfigure}
    \hfill
    \begin{subfigure}{0.4\textwidth}
        \centering
        \includegraphics[width=\textwidth]{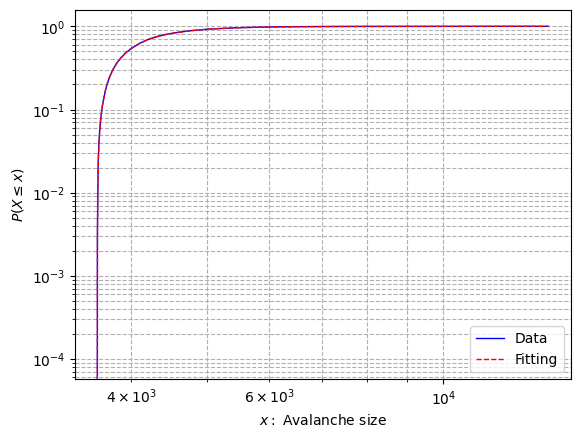}
        \caption{Model $\mathcal{G}$ of size $100\times 100$.}
        \label{fig:G100x100}
    \end{subfigure}
    \caption{Probability distribution for model $\mathcal{G}$}
    \label{fig:modelG}
\end{figure}

We consider two different evolutionary models, denoted as $\mathcal{G}$ and $\mathcal{D}$. The first model, $\mathcal{G}$, involves the evolution of the grid graph with only ``vertical'' edges to the graph grid with only ``horizontal'' edges. The evolution of $\mathcal{G}$ proceeds such that at even times, the graph consists of just ``vertical'' edges, whereas at odd times, it comprises the ``horizontal'' graph grid. 
Refer to Figure~\ref{fig:evolvingsandpilegrid3} for a visual representation of this process. 
On the other hand, the model $\mathcal{D}$ entails the regular grid graph for even times. 
However, during odd times, the neighbors of the vertex $(i,j)$ are the vertices $(i+1,j+1)$, $(i-1, j+1)$, $(i-1, j-1)$, and $(i+1, j-1)$. 
This evolution effectively tilts the edges of the regular grid graph by $45$ degrees. Refer to Figure~\ref{fig:evolvingsandpilegrid}(d) for an illustration of $\mathcal{D}$. With the proposed dynamics, we obtained a record of 10,000 avalanches with the codes provided in Appendix~\ref{codes}. Figures~\ref{fig:modelG} and~\ref{fig:modelD} depict the probability distribution of the random variable $N$ of the number of avalanches of size $x$ in models $\mathcal{G}$ and $\mathcal{D}$ respectively. We used the {\tt\color{magenta} powerlaw} python library, to fit our data to potential power law curves, see~\cite{powerlaw}. However, demonstrating that empirical data does indeed follow a power law is not an easy task. To ascertain whether our proposed evolutive sandpile models exhibit self-organized criticality in a statistical sense, we adopt the statistical framework proposed by Clauset, Shalizi, and Newman for discerning empirical power laws~\cite{ClausetPowerLaw} which is based on computing the Kolmogorov-Smirnov statistic between our empirical data and the theoretical distribution. The Kolmogorov-Smirnov statistic $\hat{\beta}$ quantifies the maximum difference between the two distributions:
\begin{equation}\label{KS}
    \hat{\beta} = \arg\min D_{\alpha}.
\end{equation}
Where $D_{\alpha}$ represents the maximum difference between the empirical distribution $F(x)$ and the power law distribution with an exponent of $\alpha$ $F_{\alpha}(x)$:
\begin{equation}\label{Dalpha}
    D_{\alpha} = \max_{x} \left|F(x)-F_{\alpha}(x)\right|.
\end{equation}

The next step proposed by Clauset, Shalizi and Newman is the \textit{log-likelihood ratio test} which calculates the value $\mathcal{R}$, representing the logarithm of the ratio between two potential distributions that fit our data. Depending on the sign of $\mathcal{R}$, we can reject or consider plausible the hypothesis that our data fits better to a distribution according to a $p$-value. Clauset, Shalizi and Newman consider it unlikely that the sign of $\mathcal{R}$ is significant if $p<0.1$. The {\tt\color{magenta} powerlaw} python library provides calculations for both the Kolmogorov-Smirnov statistic and the log-likelihood ratio test. These calculations were performed using 10,000 randomly-generated simulations. Table~\ref{tab:tabG} and Table~\ref{tab:tabD} show the values obtained from the Clausius, Shalizi and Newman framework. Based on these measurements and Figures~\ref{fig:modelG} and~\ref{fig:modelD} we can conclude that both models $\mathcal{G}$ and $\mathcal{D}$ for both sizes $50\times 50$ and $100\times 100$ indeed follow a power law.

\begin{figure}
    \centering
    \begin{subfigure}{0.4\textwidth}
        \centering
        \includegraphics[width=\textwidth]{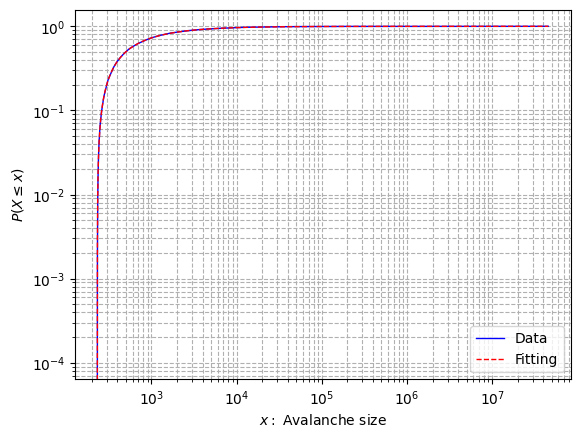}
        \caption{Model $\mathcal{D}$ of size $50\times 50$.}
        \label{fig:D50x50}
    \end{subfigure}
    \hfill
    \begin{subfigure}{0.4\textwidth}
        \centering
        \includegraphics[width=\textwidth]{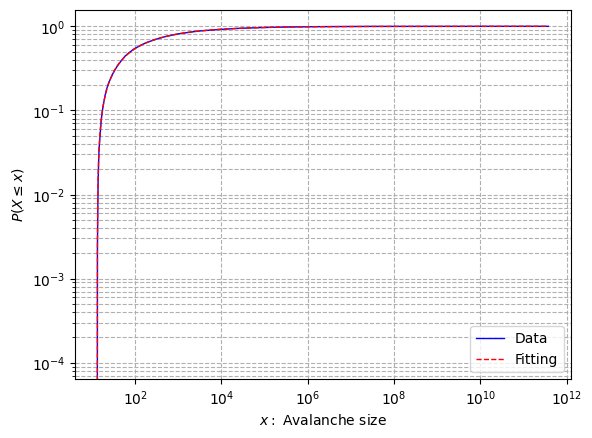}
        \caption{Model $\mathcal{D}$ of size $100\times 100$.}
        \label{fig:D100x100}
    \end{subfigure}
    \caption{Probability distribution for model $\mathcal{D}$}
    \label{fig:modelD}
\end{figure}

\begin{table}[htbp]
    \centering
    \caption{Statistic values for 10,000 simulations for model $\mathcal{G}$}
    \label{tab:tabG}
    \renewcommand{\arraystretch}{1}
    \begin{tabular}{c c c c c c}
        \arrayrulecolor{teal}
        \toprule
          & \parbox{1.5cm}{\centering Cut-off\\$x_{\text{min}}$} & \parbox{1.5cm}{\centering Exponent\\$\alpha$} & \parbox{1.5cm}{\centering Kolmogorov-Smirnov\\$\hat{\beta}$} & \parbox{1.5cm}{\centering Log-likelihood\\$\mathcal{R}$} & $p-$value\\
        \midrule
        $50\times 50$ & 853 & 8.2 & 0.006 & 0.058 & 0.95\\
        $100\times 100$ & 3622 & 8.8 & 0.006 & 0.11 & 0.91\\
        \bottomrule
    \end{tabular}
\end{table}
\begin{table}[htbp]
    \centering
    \caption{Statistic values for 10,000 simulations for model $\mathcal{D}$}
    \label{tab:tabD}
    \renewcommand{\arraystretch}{1}
    \begin{tabular}{c c c c c c}
        \arrayrulecolor{teal}
        \toprule
          & \parbox{1.5cm}{\centering Cut-off\\$x_{\text{min}}$} & \parbox{1.5cm}{\centering Exponent\\$\alpha$} & \parbox{1.5cm}{\centering Kolmogorov-Smirnov\\$\hat{\beta}$} & \parbox{1.5cm}{\centering Log-likelihood\\$\mathcal{R}$} & $p-$value\\
        \midrule
       $50\times 50$ & 231 & 1.89 & 0.004 & 0.058 & 0.95\\
       $100\times 100$ & 126 & 1.50 & 0.005 & 1.332 & 0.18\\
        \bottomrule
    \end{tabular}
\end{table}

\section{Conclusions}

Evolutive sandpile had the potential to enrich the classic sandpile models.
It is worth mentioning for instance the extension of the Schelling segregation model over evolutive graphs in \cite{henry2011emergence}.
Also, Alan Kirman has proposed that economics could be effectively modeled as a network evolving over time, wherein agents have the ability to learn from past experiences with their neighbors, allowing links between them to either weaken or strengthen~\cite{kirman1997economy}. 
In this setting, it is possible to explore classic results and extend them to evolutive sandpiles.

For example, within the algebraic perspective of the sandpiles, known as {\it sandpile group}. 
The sum of two configurations ${ c}$ and ${ d}$ is performed entry by entry.
The \textit{sandpile sum} ${ c}\oplus { d}$ of two configurations is defined as $s({ c}+{ d})$.
A configuration ${ c}$ is \textit{recurrent} if there exists a non-zero configuration ${ d}$ such that ${ c}={ c}\oplus { d}$.
Recurrent configurations play a central role in the dynamics of the Abelian sandpile model since recurrent configurations together with the sandpile sum form an Abelian group known as \textit{sandpile group} \cite[Chapter 4]{MR3889995}.
The sandpile group of $G$ is denoted by $K(G)$.
One of the interesting features of the sandpile group of connected graphs is that the order $|K(G)|$ is equal to the number $\tau(G)$ of spanning trees of the graph $G$.
For the reader interested in the mathematics of the sandpiles, we recommend the book \cite{MR3889995} or the survey \cite{AlfaroMerino}.
Is it possible to obtain such an algebraic structure on non-trivial evolutive graphs?

On the other hand, Baker and Norine introduced a \emph{dollar game} in \cite{baker2007riemann} which consists of the initial configuration assigning to each vertex $v$ an integer number of dollars. 
A vertex that has a negative number of dollars assigned to it is said to be in \emph{debt}. 
Two moves are allowed, a vertex $v$ either \emph{borrow} one dollar from each
of its neighbors or \emph{give} one dollar to each of its neighbors. 
The object of the game is to reach, through a sequence of moves, a \emph{winning configuration} that is a configuration in which no vertex is in debt. 
A sequence of moves that achieve such a configuration is called a \emph{winning strategy}.
Baker and Norine developed a graph-theoretic analog of the classical
Riemann-Roch theorem which turn out in a characterization that allows one to decide whether a winning strategy exists.
Let $g = |E(G)|- |V (G)| + 1$ and $N$ be the total number of dollars present at any stage of the game.
Therefore, (1) if $N \geq g$, then there is always a winning strategy. (2) if $N < g - 1$, then there is always an initial configuration for which no winning strategy exists.
Is it possible to obtain winning strategies in the evolutive setting?

\section*{Acknowledgement}
The authors would like to thank Lothar Dirks for his helpful comments.
The research of C.A. Alfaro is partially supported by the Sistema Nacional de Investigadores grant number 220797. The research of Ralihe R. Villagran was supported, in part, by the National Science Foundation through the DMS Award $\#$1808376 which is gratefully acknowledged.

\bibliographystyle{plain}
\bibliography{main}

\appendix
\section{Codes}\label{codes}

In this appendix, we provide the codes that were used in the manuscript.
The simulation of avalanche sizes, in addition to computing power-law distributions, was conducted using the Python programming language.

\begin{lstlisting}[language=Python, caption={Evolutive function}, label={lst:evolutive}]
import numpy
def evolutive(G: dict, c: list):
  L = list()
  t = 0
  for i in range(10000):
    avalanche_size = 0
    random_vertex = numpy.random.randint(1, len(G[t].keys()))
    c[random_vertex] += 1
    while True:
      c, k = stabilize(G[t],c)
      avalanche_size += k
      t += 1
      t %= len(G.keys())
      if k == 0: break
    L.append(avalanche_size)
  return L
\end{lstlisting}

The following demonstrates the use of the {\tt stabilize} function to stabilize sandpiles on general graphs, along with the incorporation of the auxiliary {\tt topple} method for toppling vertices.

\begin{lstlisting}[language=Python, caption={Stabilization}, label={lst:stabilize-topple}]
def topple(G: dict, c: list, v: int):
  n = len(G[v])
  grains = c[v]
  if c[v] >= n:
    for w in G[v]:
      c[w] += 1
      c[v] -= 1
  return grains >= n
  
def stabilize(G: dict, c: list):
  k = 0
  for v in range(1,len(G.keys())):
      is_topple = topple(G,c,v)
      if is_topple: k += 1
  return c, k
\end{lstlisting}





\end{document}